\documentclass[a4paper,12pt]{article}
\usepackage{amsmath,amsfonts,amssymb}
\makeatletter
\newbox\bk@bxb
\newbox\bk@bxa
\newif\if@bkcont
\newcount\bk@lcnt

\def\breakboxskip{2pt}
\def\breakboxparindent{1.8em}

\def\breakbox{\vskip\breakboxskip\relax
\setbox\bk@bxb\vbox\bgroup
\advance\linewidth -2\fboxrule
\hsize\linewidth\@parboxrestore
\parindent\breakboxparindent\relax}

\def\bk@split{%
\@tempdimb\ht\bk@bxb 
\advance\@tempdimb\dp\bk@bxb
\setbox\bk@bxa\vsplit\bk@bxb to\z@ 
\setbox\bk@bxa\vbox{\unvbox\bk@bxa}
\setbox\@tempboxa\vbox{\copy\bk@bxa\copy\bk@bxb}
\advance\@tempdimb-\ht\@tempboxa
\advance\@tempdimb-\dp\@tempboxa}

\def\bk@addfsepht{%
\setbox\bk@bxa\vbox{\vskip\fboxsep\box\bk@bxa}}

\def\bk@addskipht{%
\setbox\bk@bxa\vbox{\vskip\@tempdimb\box\bk@bxa}}

\def\bk@addfsepdp{%
\@tempdima\dp\bk@bxa
\advance\@tempdima\fboxsep
\dp\bk@bxa\@tempdima}

\def\bk@addskipdp{%
\@tempdima\dp\bk@bxa
\advance\@tempdima\@tempdimb
\dp\bk@bxa\@tempdima}

\def\bk@line{%
\hbox to \linewidth{%
\hskip-2\fboxsep\vrule \@width\fboxrule\hskip.5\fboxsep\vrule \@width\fboxrule\hskip1.5\fboxsep
\box\bk@bxa\hfil
}}%

\def\endbreakbox{\egroup
\ifhmode\par\fi{\noindent\bk@lcnt\@ne
\@bkconttrue\baselineskip\z@\lineskiplimit\z@
\lineskip\z@\vfuzz\maxdimen
\bk@split\bk@addfsepht\bk@addskipdp
\ifvoid\bk@bxb 
\def\bk@fstln{\bk@addfsepdp
\hskip-\parindent\vbox{\llap{\raisebox{-2ex}{\rule{1.5\fboxsep}{\fboxrule}\hskip.5\fboxsep}}\bk@line\llap{\rule{1.5\fboxsep}{\fboxrule}\hskip.5\fboxsep}}}

\else 
\def\bk@fstln{\vbox{\llap{\raisebox{-2ex}{\rule{1.5\fboxsep}{\fboxrule}\hskip.5\fboxsep}}\bk@line}\hfil%
\advance\bk@lcnt\@ne
\loop
\bk@split\bk@addskipdp\leavevmode
\ifvoid\bk@bxb 
\@bkcontfalse\bk@addfsepdp
\vtop{\bk@line\llap{\rule{2\fboxsep}{\fboxrule}}}%

\else 
\bk@line
\fi
\hfil\advance\bk@lcnt\@ne
\if@bkcont\repeat}%
\fi
\leavevmode\bk@fstln\par}\vskip\breakboxskip\relax}

\newcommand{\bizlie}[1]{\mathop{\,\raisebox{-.5ex}{$\widehat{\raisebox{.9ex}{\rule{2.5ex}{.07ex}}}_{#1}$}\,}}

\def\smp{\smallskip\par}
\def\un{{\bf 1}}
\def\zero{\{0\}}
\def\pf{\noindent{\bf Proof:}\ }

\def\findemo{~\leaders\hbox to 1em{\hss\  \hss}\hfill~\raisebox{.5ex}{\framebox[1ex]{}}\smp}
\def\spn{\bigskip\par\noindent}
\def\mpn{\medskip\par\noindent}
\def\smpn{\smallskip\par\noindent}
\def\normal{\mathop{\trianglelefteq}}

\def\smp{\smallskip\par}
\def\smpn{\smallskip\par\noindent}

\def\mpoint{\;\;.}
\def\mvirg{\;\;,}

\def\Res{{\rm Res}}
\def\Ind{{\rm Ind}}
\def\Inf{{\rm Inf}}
\def\Def{{\rm Def}}
\def\Iso{{\rm Iso}}

\def\Indinf{{\rm Indinf}}
\def\Defres{{\rm Defres}}

\def\Hom{{\rm Hom}}
\def\End{{\rm End}}
\def\Ext{{\rm Ext}}
\def\Inf{{\rm Inf}}
\def\Im{{\rm Im}}
\def\Out{{\rm Out}}

\def\Ker{{\rm Ker}}

\def\Id{{\rm Id}}

\def\Irr{{\rm Irr}}

\def\op{^{op}}

\def\Z{\mathbb{Z}}
\def\N{\mathbb{N}}
\def\F{\mathbb{F}}
\def\Q{\mathbb{Q}}

\newcommand{\romain}[1]{\uppercase\expandafter{\romannumeral #1}}

\newcommand{\flh}[2]{\mathop{\hbox to 12mm{\rightarrowfill}}_{\displaystyle #2}^{\displaystyle #1}\limits}
\newcommand{\sflh}[2]{\mathop{\hbox to 12mm{\rightarrowfill}}_{\scriptstyle #2}^{\scriptstyle #1}\limits}

\newcommand{\gMod}[1]{#1{\hbox{-}\mathsf{Mod}}}

\newcommand{\sur}[1]{\,\overline{\! #1}}

\def\op{^{op}}

\newcommand{\carre}[8]{\begin{array}{ccc}
#1&\mathop{\hbox to 12mm{\rightarrowfill}}^{\displaystyle{#2}}\limits&#3\\
\llap{$\displaystyle{#4}$}\left\downarrow\vbox to 6mm{}\right. & & \left\downarrow\vbox to 6mm{}\right.\rlap{$\displaystyle{#5}$}\\
#6&\mathop{\hbox to 12mm{\rightarrowfill}}_{\displaystyle #7}\limits&#8\\
\end{array}}

\newcommand{\carrem}[8]{\begin{array}{ccc}
#1&\mathop{\hbox to 12mm{\rightarrowfill}}^{\displaystyle #2}\limits&#3\\
\llap{$\displaystyle #4$}\left\uparrow\vbox to 6mm{}\right. & & \left\uparrow\vbox to 6mm{}\right.\rlap{$\displaystyle #5$}\\
#6&\mathop{\hbox to 12mm{\rightarrowfill}}_{\displaystyle #7}\limits&#8\\
\end{array}}

\newcommand{\limproj}[1]{\lim_{\displaystyle\longleftarrow\atop \scriptstyle{#1}}\limits}

\newenvironment{enonce}[1]{\pagebreak[2]\refstepcounter{subsection}\refstepcounter{prop}\smpn{{\bf \thesection.\arabic{prop}.\ \ #1:}}\begin{it} }{\end{it}\smp}
\newenvironment{enonce*}[1]{\pagebreak[2]\smpn{#1:}\begin{it} }{\end{it}\smp}
\newcommand{\result}[1]{\begin{enonce}{#1}}
\def\fresult{\end{enonce}}
\newcommand{\npar}{\smallskip\par\noindent\pagebreak[2]\refstepcounter{subsection}\refstepcounter{prop}{\bf \thesection.\arabic{prop}.\ \ }}

\newenvironment{mth}[1]{\begin{breakbox}\begin{enonce}{#1}}{\end{enonce}\end{breakbox}}
\newenvironment{mth*}[1]{\begin{breakbox}\begin{enonce*}{#1}}{\end{enonce*}\end{breakbox}}
\newenvironment{rem}[1]{\refstepcounter{subsection}\refstepcounter{prop} \mpn{{\bf \thesection.\arabic{prop}.}\ \ \bf#1:}}{\smp}

\def\dom{\backslash}
\makeatletter
\renewenvironment{enumerate}{\ifnum \@enumdepth >3 \@toodeep\else
      \advance\@enumdepth \@ne
      \edef\@enumctr{enum\romannumeral\the\@enumdepth}\list
      {\csname label\@enumctr\endcsname}{\setlength{\topsep}{1ex}\setlength{\itemsep}{0pt}\usecounter
        {\@enumctr}\def\makelabel##1{\hss\llap{##1}}}\fi}{\endlist}
\renewenvironment{itemize}{\ifnum \@itemdepth >3 \@toodeep\else \advance\@itemdepth \@ne
\edef\@itemitem{labelitem\romannumeral\the\@itemdepth}%
\list{\csname\@itemitem\endcsname}{\setlength{\topsep}{1ex}\setlength{\itemsep}{0pt}\def\makelabel##1{\hss\llap{##1}}}\fi}
{\endlist}
\makeatother
\makeatletter
\def\@sect#1#2#3#4#5#6[#7]#8{\ifnum #2>\c@secnumdepth
    \let\@svsec\@empty\else
    \refstepcounter{#1}\edef\@svsec{\csname the#1\endcsname .\hskip .5em}\fi
    \@tempskipa #5\relax
     \ifdim \@tempskipa>\z@
       \begingroup #6\relax
         \@hangfrom{\hskip #3\relax\@svsec}{\interlinepenalty \@M #8\par}%
       \endgroup
      \csname #1mark\endcsname{#7}\addcontentsline
        {toc}{#1}{\ifnum #2>\c@secnumdepth \else
                     \protect\numberline{\csname the#1\endcsname}\fi
                   #7}\else
       \def\@svsechd{#6\hskip #3\relax  
                  \@svsec #8\csname #1mark\endcsname
                     {#7}\addcontentsline
                          {toc}{#1}{\ifnum #2>\c@secnumdepth \else
                            \protect\numberline{\csname the#1\endcsname}\fi
                      #7}}\fi
    \@xsect{#5}}
\def\section{\@startsection {section}{1}{\z@}{-3.5ex plus-1ex minus
    -.2ex}{2.3ex plus.2ex}{\reset@font\Large\bf}}  

\makeatother
\renewenvironment{equation}{\refstepcounter{subsection}\refstepcounter{prop}$$}{\leqno{\bf (\theprop)}$$}

\def\mar[#1]{\ar@{-}[#1]|-{\object@{<}}}
\def\marb[#1]{\ar@{-}[#1]|{\object+{  }}}

\usepackage[all]{xy}
\def\endpf{\findemo}
\def\CC{\mathcal{C}}
\def\CE{\mathcal{E}}
\def\CF{\mathcal{F}}
\def\CG{\mathcal{G}}
\def\CO{\mathcal{O}}
\def\CR{\mathcal{R}}
\def\CS{\mathcal{S}}
\def\CT{\mathcal{T}}
\begin{document}
\centerline{\Large \bf A functorial presentation}\vspace{2ex}
\centerline{\Large \bf of units of Burnside rings}\vspace{3ex}\par
\centerline{\large\bf Serge Bouc}\vspace{3ex}
{\footnotesize
{\bf Abstract:} Let $B^\times$ be the biset functor over $\F_2$ sending a finite group~$G$ to the group $B^\times(G)$ of units of its Burnside ring $B(G)$, and let $\widehat{B^\times}$ be its dual functor. The main theorem of this paper gives a characterization of the cokernel of the natural injection from $B^\times$ in the dual Burnside functor $\widehat{\F_2B}$, or equivalently, an explicit set of generators $\CG_\CS$ of the kernel $L$ of the natural surjection $\F_2B\to \widehat{B^\times}$. This yields a two terms projective resolution of $\widehat{B^\times}$, leading to some information on the extension functors $\Ext^1(-,B^\times)$. For a finite group~$G$, this also allows for a description of $B^\times(G)$ as a limit of groups $B^\times(T/S)$ over sections $(T,S)$ of $G$ such that $T/S$ is cyclic of odd prime order, Klein four, dihedral of order 8, or a Roquette 2-group. Another consequence is that the biset functor $B^\times$ is not finitely generated, and that its dual $\widehat{B^\times}$ is finitely generated, but not finitely presented. The last result of the paper shows in addition that $\CG_\CS$ is a minimal set of generators of~$L$, and it follows that the lattice of subfunctors of $L$ is uncountable.}\vspace{1ex}\par
\par

\section{Introduction}
The Burnside ring $B(G)$ is a fundamental invariant attached to a finite group~$G$. It is a commutative (unital) ring, and most of its structural properties have been described several decades ago, e.g. the prime spectrum by Dress (\cite{dressresoluble}), or the primitive idempotents of the algebra $\Q B(G)$ by Gluck (\cite{gluck}) and Yoshida (\cite{yoshidaidemp}). \par
An important missing item in this list is the group of multiplicative units $B(G)^\times$. More precisely, it follows from Burnside's theorem that $B(G)^\times$ is a finite elementary abelian 2-group, but the rank of this group is known only under additional assumptions on $G$, by the work of many different people (\cite{matsuda}, \cite{matsudamiyata}, \cite{yoshidaunit}, \cite{yalcin}, \cite{burnsideunits}, \cite{barsotti}). A very efficient algorithm for computing this rank has also been obtained in \cite{boltje-pfeiffer}. However, no general formula for the rank of $B(G)^\times$ is known so far. Let us recall that, according to an observation of tom Dieck (\cite{tomdieckgroups}, Proposition 1.5.1) based on a theorem of Dress (\cite{dressresoluble}), Feit-Thompson's theorem (\cite{feit-thompson}) is equivalent to the assertion that $B^\times(G)$ has order 2 if $|G|$ is odd.\par
It was observed in~\cite{burnsideunits} that the assignment $G\mapsto B(G)^\times$ is a biset functor $B^\times$ with values in $\F_2$-vector spaces, and that $B^\times$ embeds in the $\F_2$-dual functor $\Hom_\Z(B,\F_2)$ of the Burnside  functor. In the present paper, we give (Theorem~\ref{image}) a characterization of the image of this embedding $\imath:B^\times\to\Hom_\Z(B,\F_2)$, or equivalently, we describe generators for the kernel $L$ of the natural transposed surjection $\jmath:\F_2B\to\Hom_\Z(B^\times,\F_2)$. This yields a functorial presentation
$$\xymatrix{0\ar[r]&L\ar[r]&\F_2B\ar[r]^-\jmath&\Hom_\Z(B^\times,\F_2)\ar[r]&0}$$
of the (dual) functor $\Hom_\Z(B^\times,\F_2)$. Unfortunately, the generators we obtain for $L$ are not linearly independent, so at least at this stage, the previous exact sequence doesn't provide any obvious formula for the $\F_2$-dimension of the evaluations of the functor $B^\times$.\par
It still admits some interesting consequences, that we develop in the three last sections. Section~\ref{ext1} is concerned with the extension groups $\Ext^1(M,B^\times)$ for some biset functors $M$ over $\F_2$. In Section~\ref{char}, we give a sectional characterization of $B^\times(G)$ as a limit of groups $B^\times (T/S)$, where $(T,S)$ runs through sections $(T,S)$ of $G$ such that $T/S$ is cyclic of odd prime order, Klein four, dihedral of order 8, or a Roquette 2-group. Finally, in Section 6, it is shown that the biset functor $B^\times$ is not finitely generated, and that its dual is finitely generated (by a single element), but not finitely presented. We also show that the generating set we obtain for the above functor $L$ is minimal, and that the lattice of subfunctors of $L$ is uncountable. These last results may be seen as another sign of the difficulty of determining $B^\times(G)$ for an arbitrary finite group $G$, at least with our present methods.
\section{Review of Burnside rings and biset functors}
\npar {\bf Burnside rings.} We quickly recall first some basic definitions on Burnside rings. Missing details can be found in~\cite{handbook}.\par
Let $G$ be a finite group. The Burnside ring $B(G)$ of $G$ is the Grothendieck ring of the category of finite (left) $G$-sets, for relations given by disjoint union decompositions. In other words $B(G)$ is the quotient of the free abelian group on the set of isomorphism classes of finite $G$-sets by the subgroup generated by all the elements of the form $[X\sqcup Y]-[X]-[Y]$, where $[X]$ denote the isomorphism class of the $G$-set $X$. \par
The multiplication in $B(G)$ is induced by the cartesian product of $G$-sets. In other words $[X][Y]=[X\times Y]$ for any two finite $G$-set $X$ and $Y$. The identity element is the (isomorphism class of a) $G$-set of cardinality 1. \par
The additive group $B(G)$ is free with basis the set of isomorphism classes of transitive $G$-sets. Each transitive $G$-set is isomorphic to a $G$-set of the form $G/H$, where $H$ is a subgroup of $G$. The transitive $G$ sets $G/H$ and $G/K$ are isomorphic if and only if the subgroups $H$ and $K$ of $G$ are conjugate. For sake of simplicity, we denote by $G/H$ instead of $[G/H]$ the image in $B(G)$ of the transitive $G$-set $G/H$. With this abuse of notation, the abelian group $B(G)$ has a basis consisting of the elements $G/H$, where $H$ runs through a set of representatives of conjugacy classes of subgroups of $G$.\par
The commutative $\Q$-algebra $\Q B(G)=\Q\otimes_\Z B(G)$ is split semisimple. Its primitive idempotents (see \cite{gluck} or \cite{yoshidaidemp}) are indexed by the subgroup of $G$, up to conjugation. The idempotent $e_H^G$ indexed by $H\leq G$ is equal to
$$e_H^G=\frac{1}{|N_G(H)|}\sum_{L\leq H}|L|\mu(L,H)\,G/L\mvirg$$
where $\mu$ is the M\"obius function of the poset of subgroups of $G$. Its defining property is that
$$\forall K\leq G,\;|(e_H^G)^K|=\left\{\begin{array}{ll}1&\hbox{if}\;K=_GH,\\0&\hbox{otherwise.}\end{array}\right.$$
Here the map $\alpha\in\Q B(G)\mapsto |\alpha^K|\in\Q$ is the linear form sending (the class of a finite) $G$-set $X$ to the cardinality $|X^K|$ of the set $X^K|$ of $K$-fixed points in $X$.
\npar {\bf Biset functors.} (see \cite{bisetfunctorsMSC} for details) For finite groups $G$ and $H$, an $(H,G)$-biset is a $H\times G\op$-set, i.e. a set $U$ endowed with a left $H$-action and a right $G$-action which commute. The Burnside group $B(H\times G\op)$ of (finite) $(H,G)$-bisets is denoted by $B(H,G)$. \par
Let $\mathcal{C}$ denote the following category:
\begin{itemize}
\item The objects of $\mathcal{C}$ are all the finite groups.
\item For finite groups $G$ and $H$, the set of morphisms $\Hom_\mathcal{C}(G,H)$ is equal to $B(H,G)$.
\item The composition of morphisms in $\mathcal{C}$ is defined by bilinearity from the standard ``tensor product'' (also called {\em composition}) of bisets: for finite groups $G$, $H$ and $K$, for a $(K,H)$-biset $V$ and an $(H,G)$-biset $U$, set 
$$V\times_HU=(V\times U)/H\mvirg$$
where $H$ acts on the right on $(V\times U)$ by $(v,u)h=(vh,h^{-1}u)$. Then $V\times_HU$ is a $(K,G)$-biset in the obvious way.
\item The identity element of $G$ is (the class of) the set $G$, viewed as a $(G,G)$-biset by left and right multiplication.
\end{itemize}
A {\em biset functor} is an additive functor from $\mathcal{C}$ to the category of all $\Z$-modules. Biset functors, together with natural transformations between them, form an abelian category $\mathcal{F}$. More generally, for a commutative (unital) ring $k$, one can consider the $k$-linearization $k\mathcal{C}$ of $\mathcal{C}$, i.e. the category with the same objects, and such that 
$$\Hom_{k\mathcal{C}}(G,H)=k\otimes_\Z\Hom_\mathcal{C}(G,H)=kB(H,G)\mvirg$$
the composition in $k\CC$ being defined as the $k$-bilinear extension of the composition in $\mathcal{C}$. A {\em biset functor over $k$} is a $k$-linear functor from $k\mathcal{C}$ to the category $\gMod{k}$ of all $k$-modules. These functors, together with natural transformations between them, form an abelian category $\mathcal{F}_k$. \par
There is a one to one parametrization $(H,V)\mapsto S_{H,V}$ of simple bisets functors over $k$ (up to isomorphism of functors) by pairs $(H,V)$ of a finite group $H$ and a simple $k\Out(H)$-module $V$ (up to isomorphism of such pairs). If $G$ is a finite group, and if $S_{H,V}(G)\neq\zero$, then $H$ is isomorphic to a subquotient of $G$.\par
Let $F$ be a biset functor, and let $G$, $H$ be finite groups. We define {\em the opposite biset} $U\op$ as the $(G,H)$-biset $U$, where the action of $(g,h)\in G\times H$ on $u\in U$ is defined by
$$g\cdot u\cdot h=h^{-1}ug^{-1}\mpoint$$
This definition extends uniquely to a $k$-linear map $\alpha\mapsto \alpha\op$ from $kB(H,G)$ to $kB(G,H)$.\par
For $\alpha\in kB(H,G)$ and $m\in F(G)$, we denote by $\alpha(m)$ or simply $\alpha\, m$ the image of $m$ by the map $F(\alpha)$. In particular, when $U$ is a finite $(H,G)$-biset, we set $U\,m=F\big([U]\big)(m)$. \par
When $F$ is a biset functor over $k$, and $M$ is a $k$-module, the $M$-dual of~$F$ is the functor $\Hom_k(F,M)$ defined by
$$\left\{\begin{array}{l}\forall G,\;\Hom_k(F,M)(G)=\Hom_k\big(F(G),M\big)\\\forall G,H,\,\forall \alpha\in kB(H,G),\;\Hom_k(F,M)(\alpha)={^tF(\alpha\op)}\mpoint\end{array}\right.$$
When $p$ is a prime number, let $k\CC_p$ denote the full subcategory of $k\CC$ consisting of finite $p$-groups. A $k$-linear functor from $k\CC_p$ to $\gMod{k}$ is called a {\em $p$-biset functor} (over $k$). The abelian category of $p$-biset functors over $k$ is denoted by $\CF_{p,k}$.\par
More generally, if $\mathcal{T}$ is a class of finite groups closed under taking subquotients, one can consider the full subcategory $k\mathcal{T}$ of $k\CC$ consisting of groups in~$\mathcal{T}$. The abelian category of $k$-linear functors from $k\mathcal{T}$ to $\gMod{k}$ is denoted by $\CF_{\CT,k}$.\par
Let $F$ be an object of $\CF_{T,k}$. An {\em element} of $F$ is a pair $(R,v)$, where $R\in\CT$ and $v\in F(R)$. For a set $\mathcal{G}$ of elements of $F$, the subfunctor $\langle \mathcal{G}\rangle$ of $F$ {\em generated} by $\mathcal{G}$ is defined as the intersection of all subfunctors $L$ of $F$ such that $v\in L(R)$ for all $(R,v)\in \mathcal{G}$. Its evaluation at a group $H\in\CT$ is given by
$$\langle \mathcal{G}\rangle(H)=\sum_{\substack{(R,v)\in \mathcal{G}\\\alpha\in kB(H,R)}}F(\alpha)(v)\mpoint$$
The functor $F$ is {\em finitely generated} if there exists a finite set $\mathcal{G}$ of elements of~$F$ such that $\langle \mathcal{G}\rangle=F$. The functor $F$ is {\em finitely presented} if there exists an exact sequence $Q\to P\to F\to 0$ in $\CF_{\CT,k}$, where  $P$ and $Q$ are finitely generated projective functors.\par
Let $F$ be an object of $\CF_{T,k}$, and $H\in\CT$. The {\em residue} (or {\em Brauer quotient}) $\sur{F}(H)$ of $F$ at $H$ is the $k$-module defined by
$$\sur{F}(H)=F(H)\big/\sum_{\substack{T\in\CT\\|T|<|H|\\\alpha\in kB(H,T)}}F(\alpha)\big(F(T)\big)\mpoint$$
\npar{\bf Elementary bisets.} Bisets of one of the following forms are called {\em elementary bisets}:
\begin{itemize}
\item Let $H$ be a subgroup of a finite group $G$. The set $G$, viewed as a $(G,H)$-biset (by multiplication), is called {\em induction} (from $H$ to $G$), and denoted by $\Ind_H^G$. The set $G$, viewed as an $(H,G)$-biset (by multiplication), is called {\em restriction} (from $G$ to $H$), and denoted by $\Res_H^G$.
\item Let $N$ be a normal subgroup of a finite group $G$. The set $G/N$, viewed as a $(G,G/N)$-biset (by multiplication on the right, and projection followed by multiplication on the left), is called {\em inflation} (from $G/N$ to~$G$) and denoted by $\Inf_{G/N}^G$. The set $G/N$, viewed as a $(G/N,G)$-biset (by multiplication on the left, and projection followed by multiplication on the right), is called {\em deflation} (from $G$ to $G/N$), and denoted by $\Def_{G/N}^G$.
\item Let $\varphi:G\to G'$ be an isomorphism of groups. The set $G'$, viewed as a $(G',G)$-biset (by multiplication on the left, and multiplication by the image under $\varphi$ on the right), is called {\em transport by isomorphism} (by $\varphi$), and denoted by $\Iso(\varphi)$ or $\Iso_{G}^{G'}$ if $\varphi$ is clear from the context.
\end{itemize}
In addition to these elementary bisets, it is convenient to distinguish the following ones, when $(T,S)$ is a section of a finite group $G$, i.e. a pair of subgroups of $G$ with $S\normal T$:
\begin{itemize}
\item The set $G/S$, viewed as a $(G,T/S)$-biset, is called {\em induction-inflation} (from $T/S$ to $G$), and denoted by $\Indinf_{T/S}^G$. It is isomorphic to the composition $\Ind_T^G\times_T\Inf_{T/S}^T$. 
\item The set $S\dom G$, viewed as a $(T/S,G)$-biset, is called {\em deflation-restriction} (from $G$ to $T/S$), and denoted by $\Defres_{T/S}^G$. It is isomorphic to the composition $\Def_{T/S}^T\times_T\Res_T^G$.
\end{itemize}
For finite groups $G$ and $H$, any transitive $(H,G)$-biset is isomorphic to a biset of the form $(H\times G)/X$, where $X$ is a subgroup of $H\times G$, and the biset structure is given by $h\cdot (a,b)X\cdot g=(ha,g^{-1}b)X$ for $h,a$ in $H$ and $g,b$ in $G$. We set
\begin{align*}
p_1(X)&=\{h\in H\mid\exists g\in G,\;(h,g)\in X\}\\
k_1(X)&=\{h\in H\mid(h,1)\in X\}\\
p_2(X)&=\{g\in G\mid\exists h\in H,\;(h,g)\in X\}\\
k_2(X)&=\{g\in G\mid(1,g)\in X\}\mpoint\\
\end{align*}
With this notation, we have $k_1(X)\normal p_1(X)$ and $k_2(X)\normal p_2(X)$, and there is a canonical group isomorphism $f:p_2(X)/k_2(X)\to p_1(X)/k_1(X)$ sending $gk_2(X)$ to $hk_1(X)$ if $(h,g)\in X$. Moreover (see \cite{bisetfunctorsMSC}, Lemma 2.3.26), there is an isomorphism of $(H,G)$-bisets
\begin{equation}\label{decomposition}(H\times G)/X\cong \Indinf_{p_1(X)/k_1(X)}^H\,\Iso(f)\,\Defres_{p_2(X)/k_2(X)}^G
\end{equation}
where the concatenation on the right hand side denotes the composition of bisets. In other words, any transitive biset is isomorphic to a composition of elementary bisets.\par
It follows that if $\CT$ is a class of finite groups closed under taking subquotients, if $F\in\CF_{\CT,k}$ and $H\in\CT$, then
\begin{align*}\sur{F}(H)&=F(H)\big/\sum_{\substack{A\normal B\leq H\\(B,A)\neq (H,\un)}}\Indinf_{B/A}^HF(B/A)\\
&=F(H)\big/\Big(\sum_{\substack{B<H\\B\,\mathrm{maximal}}}\Ind_B^HF(B)+\sum_{\substack{\un<A\normal H\\A\,\mathrm{minimal}}}\Inf_{H/A}^HF(H/A)\Big)\mpoint
\end{align*}
\npar {\bf Faithful elements.} (see \cite{bisetfunctorsMSC}, Section 6.3) Let $F$ be a biset functor over a commutative ring $k$. For a finite group $G$, let
$$\partial F(G)=\{u\in F(G)\mid \forall \un\neq N\normal G,\; \Def_{G/N}^Gu=0\}\mpoint$$
The $k$-submodule $\partial F(G)$ is called the submodule of {\em faithful elements} of $F(G)$. If is always a direct summand of $F(G)$. More precisely, the element
$$f_\un^G=\sum_{N\normal G}\mu_{\normal G}(\un,N)\,[(G\times G)/L]$$
of $kB(G,G)$ is an idempotent endomorphism of $G$ in the category $k\CC$, and one can show that $\partial F(G)=F(f_\un^G)\big(F(G)\big)$.
\pagebreak[3]
\npar {\bf Genetic bases of $p$-groups.} (\cite{bisetfunctorsMSC} Definition 6.4.3, Lemma 9.5.2 and Theorem 9.6.1) Let $p$ be a prime number, and $P$ be a finite $p$-group. For a subgroup $Q$ of $P$, let $Z_P(Q)\geq Q$ denote the subgroup of $N_P(Q)$ defined by
$$Z_P(Q)/Q=Z\big(N_P(Q)/Q\big)\mpoint$$
The subgroup $Q$ is a {\em genetic} subgroup of $P$ if the two following properties hold:
\begin{itemize}
\item The group $N_P(Q)/Q$ is a Roquette $p$-group, i.e. it has {\em normal $p$-rank~1}. Recall that the Roquette $p$-groups of order $p^n$ are the cyclic $p$-groups $C_{p^n}$ ($n\geq 0$), and in addition when $p=2$, the generalized quaternion groups $Q_{2^n}$ ($n\geq 3$), the dihedral groups $D_{2^n}$ ($n\geq 4$) and the semidihedral groups $SD_{2^n}$ ($n\geq 4$).
\item For any $x\in G$, the intersection ${^xQ}\cap Z_P(Q)$ is contained in $Q$ if and only if $^xQ=Q$.
\end{itemize} 
For two genetic subgroups $Q$ and $R$ of $P$, write
$$Q\bizlie{P}R \Leftrightarrow \exists x\in P \;\hbox{such that}\; ^xQ\cap Z_P(R)\leq R\; \hbox{and}\; R^x\cap Z_P(Q)\leq Q\mpoint$$ 
One can show that this defines an equivalence relation on the set of genetic subgroups of $P$. {\em A genetic basis} of $P$ is a set of representatives of equivalence classes of genetic subgroups of $P$ for the relation $\bizlie{P}$.
\npar {\bf Rational $p$-biset functors.} (\cite{bisetfunctorsMSC} Theorem~10.1.1, Definition 10.1.3 and Theorem 10.1.5) Let $p$ be a prime number, and $F$ be a $p$-biset functor over a commutative ring $k$. If $P$ is a finite $p$-group and $\mathcal{B}$ is a genetic basis of $P$, the map
$$\mathcal{I}_\mathcal{B}:\bigoplus_{Q\in\mathcal{B}}\Indinf_{N_P(Q)/Q}^P:\bigoplus_{Q\in\mathcal{B}}\partial F\big(N_P(Q)/Q\big)\to F(P)$$
is split injective, with left inverse
$$\mathcal{D}_\mathcal{B}:\bigoplus_{Q\in\mathcal{B}}f_\un^{N_P(Q)/Q}\circ\Defres_{N_P(Q)/Q}^P:F(P)\to\bigoplus_{Q\in\mathcal{B}}\partial F\big(N_P(Q)/Q\big)\mpoint$$
The functor $F$ is called {\em rational} if for any finite $p$-group $P$, the map $\mathcal{I}_\mathcal{B}$ is an isomorphism for {\em some} - equivalently for {\em any} - genetic basis $\mathcal{B}$ of $P$. Rational $p$-biset functors form a Serre subcategory of $\CF_{p,k}$. Moreover, for any $k$-module~$M$, the $M$-dual $\Hom_k(F,M)$ of a rational $p$-biset functor $F$ is rational.
\section{The main theorem}
\begin{mth}{Notation} Let $\mathcal{R}$ denote the class of finite groups which are isomorphic to one of the following groups:
\begin{itemize}
\item A cyclic group $C_p$ of odd prime order $p$.
\item A cyclic group $C_4$ of order 4.
\item An elementary abelian group $(C_2)^2$ of order 4.
\item A dihedral 2-group $D_{2^n}$ of order $2^n\geq 8$.
\item A semidihedral 2-group $SD_{2^n}$ of order $2^n\geq 16$.
\end{itemize}
For each group $R$ in $\mathcal{R}$, let $\varepsilon_R$ denote the element of $B(R)$ defined by:
\begin{itemize}
\item	$\varepsilon_R=R/\un-R/R$ if $R\cong C_p$.
\item $\varepsilon_R=R/\un-R/S$ if $R\cong C_4$, where $S$ is the unique subgroup of order 2 of $R$.
\item $\varepsilon_R=R/\un-(R/A+R/B+R/C)+2R/R$ if $R\cong (C_2)^2$, where $A,B,C$ are the the three subgroups of order 2 of $R$.
\item $\varepsilon_R=(R/I-R/IZ)-(R/J-R/JZ)$, if $R\cong D_{2^n}$, where $I$ and $J$ are non conjugate non central subgroups of order 2 of $R$, and $Z$ is the center of~$R$.
\item $\varepsilon_R=R/I-R/IZ$ if $R\cong SD_{2^n}$, where $I$ is a non central subgroup of order 2 of $R$, and $Z$ is the center of~$R$.
\end{itemize}
Let moreover $\sur{\varepsilon}_R$ denote the image of $\varepsilon_R$ in $\F_2B(G)$.
\end{mth}
\begin{rem}{Remark} \label{invariant} 
When $R$ is cyclic, elementary abelian of order 4, or dihedral, the elements $\varepsilon_R\in B(R)$ already appear in \cite{boya} (Notation 3.6) and \cite{dadegroup} (Corollary 6.5 and Notation 6.9, where $\varepsilon_{D_{2^n}}$ is denoted $\delta_n$). Also observe that $\varepsilon=f_\un^R\,R/\un\in B(R)$ if $R$ is neither dihedral nor semidihedral. Moreover $\varepsilon_R=f_\un^R\,(R/I-R/J)$ if $R\cong D_{2^n}$, and $\varepsilon_R=f_\un^R\,R/I$ when $R\cong SD_{2^n}$. Hence $f_\un^R\,\varepsilon_R=\varepsilon_R$ for any $R\in\mathcal{R}$.
\end{rem}
\pagebreak[3]
\begin{mth}{Notation} Let $\F_{2,+}$ denote the additive group of $\F_2$, and 
$$s\in\{\pm1\}\mapsto s_+\in \F_{2,+}$$
be the group isomorphism is defined by $(-1)_+=1_{\F_2}$ and $(+1)_+=0_{\F_2}$. Let $t\in \F_{2,+}\mapsto t_\times\in\{\pm1\}$ denote the inverse group isomorphism.
\end{mth}
\pagebreak[3]
Recall that for a finite group $G$, any unit element $u\in B^\times(G)$ in the Burnside ring of $G$ defines a linear form $\imath_G(u)$ on $B(G)$, with values in $\F_2$, by
$$\imath_G(u)(G/H)=|u^H|_+\mpoint$$
When $G$ runs through finite groups, these maps $\imath_G$ form an injective morphism of biset functors $\imath:B^\times\to \Hom_\Z(B,\F_2)$ (\cite{burnsideunits}, Proposition 7.2). \par
Moreover, by a theorem of Yoshida (\cite{yoshidaunit}, Proposition 6.5), a linear form $\varphi: B(G)\to\F_2$ lies in the image of $\imath_G$ if and only if for any subgroup $H$ of~$G$, the map 
$$\widetilde{\varphi}_H:x\in N_G(H)/H\mapsto \varphi\big(G/H\langle x\rangle\big)-\varphi(G/H)$$
is a group homomorphism from $N_G(H)/H$ to $\F_{2,+}$.
\begin{mth}{Theorem} \label{image}Let $G$ be a finite group, and $\varphi\in\Hom_\Z\big(B(G),\F_2\big)$. The following assertions are equivalent:
\begin{enumerate}
\item $\varphi$ lies in the image of the natural map $\imath_G:B^\times(G)\to \Hom_\Z\big(B(G),\F_2\big)$.
\item $(\Defres_{T/S}^G\varphi)(\varepsilon_{T/S})=0_{\F_2}$ for every section $(T,S)$ of $G$ such that $T/S\in\nolinebreak\mathcal{R}$.
\end{enumerate}
\end{mth}
\pf \fbox{$1\Rightarrow 2$} We check that if $\varphi\in\Im\,\imath_G$, then $(\Defres_{T/S}^G\varphi)(\varepsilon_{T/S})=0_{\F_2}$ for each section $(T,S)$ of $G$ such that $T/S\in\mathcal{R}$. Since $\Defres_{T/S}^G\varphi\in\Im\,\imath_{T/S}$ as $\imath$ is a morphism of biset functors, we can assume that $G\in\mathcal{R}$, and check that $\varphi(\varepsilon_G)=0_{\F_2}$ when $\varphi=\imath_G(u)$ for some $u\in B^\times(G)$. In this case $\varphi(G/H)=|u^H|_+$ for any subgroup $H$ of $G$.\par
But since $f_\un^G\,\varepsilon_G=\varepsilon_G$ by Remark~\ref{invariant}, and since $f_\un^G=(f_\un^G)\op$, we have
\begin{align*}
\varphi(\varepsilon_G)&=i_G(u)(f_\un^G\,\varepsilon_G)\\
&=\big(f_\un^Gi_G(u)\big)(\varepsilon_G)\\
&=i_G(f_\un^G\,u)(\varepsilon_G)\mpoint
\end{align*}
Moreover $f_\un^G\,u\in\partial B^\times(G)$ by definition. But $\partial B^\times(G)=\zero$ if $|Z(G)|>2$ by Lemma~6.8 of~\cite{burnsideunits}, and $\partial B^\times(G)$ also vanishes if $G$ is generalized quaternion or semidihedral, by Corollary~6.10 there. It follows that $\partial B^\times(G)=\zero$ if $G\in \CR$, unless perhaps if $G\cong D_{2^n}$ for some $n\geq 3$. In particular $\varphi(\varepsilon_G)=0_{\F_2}$ if $G\in\CR$, unless $G$ is a dihedral 2-group of order at least 8.\par
Now by Corollary~6.12 of~\cite{burnsideunits} (and its proof), if $G\cong D_{2^n}$ for $n\geq 3$, the group $\partial B^\times(G)$ has order 2, generated by the element
$$\upsilon_G=1-2(e_I^G+e_J^G)=G/G+G/\un-(G/I+G/J)$$
of $B(G)$. Moreover $\varepsilon_G=f_\un^G\,(G/I-G/J)$ by Remark~\ref{invariant}. Thus
\begin{align*}
\varphi(\varepsilon_G)&=i_G(u)\big(f_\un^G\,(G/I-G/J)\big)\\
&=\big(f_\un^Gi_G(u)\big)(G/I-G/J)\\
&=i_G(f_\un^G\,u)(G/I-G/J)\\
&=|(f_\un^G\,u)^I|_+-|(f_\un^G\,u)^J|_+\mpoint
\end{align*}
We can moreover assume that $f_\un^G\,u$ is the generator $\upsilon_G$ of $\partial B^\times(G)$, i.e. that $f_\un^G\,u=1-2(e_I^G+e_J^G)$. Then clearly $|(f_\un^G\,u)^I|=|(f_\un^G\,u)^J|$, so $\varphi(\varepsilon_G)=0_{\F_2}$ also in this case, as was to be shown.

\spn
\fbox{$2\Rightarrow 1$} We prove the converse by induction on the order of $G$: we consider $\varphi\in\Hom_\Z\big(B(G),\F_2\big)$ such that $(\Defres_{T/S}^G\varphi)(\varepsilon_{T/S})=0_{\F_2}$ for every section $(T,S)$ of $G$ such that $T/S\in\mathcal{R}$. Assuming the result holds for groups of order smaller than $|G|$, by transitivity of deflation-restrictions, we can assume that $\Defres_{T/S}^G\varphi$ lies in the image of $B^\times(T/S)$ for any section $(T,S)$ such that $|T/S|<|G|$, in other words any section $(T,S)$ different from $(G,1)$. By Yoshida' theorem above, it follows that the map $\widetilde{\varphi}_H$ is a group homomorphism for any $H\neq \un$. So proving that $\varphi$ lies in the image of $\imath_G$ amounts to proving that the map $\widetilde{\varphi}_\un:G\to \F_2$ is a group homomorphism.\vspace{.5ex}\par
Equivalently, all we have to show is that the map 
$$\psi:x\in G\mapsto \widetilde{\varphi}_\un(x)_\times\in\{\pm1\}$$
is a group homomorphism. Now $\psi$ is a central function on $G$, with values in $\{\pm1\}$, and moreover $\psi(1)=1$ since $\widetilde{\varphi}_1(1)=0_{\F_2}$. Hence all we have to show is that $\psi$ is a generalized character of $G$. Indeed in this case, there exist integers $n_\chi\in\Z$, indexed by the irreducible complex characters of $G$ such that $\psi=\mathop{\sum}_{\chi\in\Irr(G)}\limits n_\chi\chi$. Moreover
$$\sum_{\chi\in\Irr(G)}n_\chi^2=\frac{1}{|G|}\sum_{x\in G}|\psi(x)|^2=1\mvirg$$
so $n_\chi=0$, except for a single irreducible character $\chi$ for which $n_\chi=\pm1$. Since $\psi(1)=1$, it follows that $\psi=\chi$, so $\psi$ is a character of $G$, of degree 1, hence a group homomorphism from $G$ to $\{\pm1\}$.\vspace{.5ex}\par
Now by Brauer's induction theorem, the map $\psi$ is a generalized character of $G$ if and only if its restriction to any Brauer elementary subgroup $H$ of $G$ is a generalized character of $H$. But for $x\in H$
\begin{align*}
(\Res_H^G\psi)(x)=\widetilde{\varphi}_1(x)_\times&=\varphi(G/\langle x\rangle)_\times/\varphi(G/\un)_\times\\
&=\varphi\Big(\Ind_H^G\big(H/\langle x\rangle-H/\un\big)\Big)_\times\\
&=\Big((\Res_H^G\varphi)\big(H/\langle x\rangle-H/\un\big)\Big)_\times\mpoint
\end{align*}
It follows that if $H$ is a proper subgroup of $G$, since $\Res_H^G\varphi\in\Im\,\imath_H$, the map $\Res_H^G\psi$ is a generalized character of $H$. In other words we can assume that $G$ itself is Brauer elementary.\vspace{.5ex}\par
Then there exists a prime number $p$ such that $G=C\times P$, where $P$ is a $p$-group, and $C$ is a $p'$-cyclic group. In particular $G$ is nilpotent, so $G=Q\times R$, where $Q$ is a 2-group and $R$ is nilpotent of odd order. Then any subgroup $T$ of $G$ is equal to $A\times B$, for some subgroup $A$ of $Q$ and some subgroup $B$ of~$R$. If $B$ is non trivial, there exists a normal subgroup $C$ of $B$ of (odd) prime index $l$. Set $S=A\times C$. Then $(T,S)$ is a section of $G$, and $T/S\cong C_l$. Thus
$$(\Defres_{T/S}^G\varphi)(\varepsilon_{T/S})=0_{\F_2}=\varphi\big(\Indinf_{T/S}^G(T/S-T/T)\big)=\varphi(G/S)-\varphi(G/T)\mpoint$$
It follows by induction that $\varphi\big(G/(A\times B))=\varphi\big(G/(A\times \un)\big)$ for any subgroup $A$ of $Q$ and any subgroup $B$ of $R$. In particular, for $x\in Q$ and $y\in R$, the subgroup $\langle (x,y)\rangle$ of $G=Q\times R$ is equal to $\langle x\rangle\times \langle y\rangle$, and
\begin{align*}
\psi(x,y)&=\varphi\big(G/(\langle x\rangle\times \langle y\rangle)\big)_\times/\varphi(G/\un)_\times\\
&= \varphi\big(G/(\langle x\rangle\times \un\rangle)\big)_\times/\varphi(G/\un)_\times\\
&=\psi(x,1)\mpoint
\end{align*}
In other words it is enough to prove that the restriction of $\psi$ to $Q\times \un$ is a generalized character. Equivalently, all we have to do is consider the case where $G=Q\times \un$, i.e. we can assume that $G$ is a 2-group.\vspace{.5ex}\par
We consider the linearization morphism $B(G)\to R_\Q(G)$, which is surjective by the Ritter-Segal theorem (\cite{ritter}, \cite{segal}, \cite{rittersegal}), and hence fits in a short exact sequence
$$\xymatrix{
0\ar[r]&K(G)\ar[r]&B(G)\ar[r]&R_\Q(G)\ar[r]&0
}\mpoint$$
By Corollary~6.16 of~\cite{dadegroup}, the kernel $K(G)$ of the linearization morphism is linearly generated by the elements of the form $\Indinf_{T/S}^G\varepsilon_{T/S}$, where $(T,S)$ is a section of $G$ such that $T/S$ is elementary abelian of order 4 or dihedral of order at least 8. Since 
$$\varphi(\Indinf_{T/S}^G\varepsilon_{T/S})=(\Defres_{T/S}^G\varphi)(\varepsilon_{T/S})=0_{\F_2}\mvirg$$
the linear form $\varphi$ vanishes on $K(G)$, so we can consider $\varphi$ as an element of $\Hom_\Z\big(R_\Q(G),\F_2\big)$. \par
Now the functor $R_\Q$ is a {\em rational 2-biset functor} (see Definition~10.1.3 and Proposition~9.6.12 of~\cite{bisetfunctorsMSC}), so its $\F_2$-dual $\Hom_\Z(R_\Q,\F_2)$ is also rational. It follows that if $\mathcal{G}$ is a genetic basis of $G$, we have
\begin{equation}\label{eq5}\varphi=\sum_{S\in\mathcal{G}}\Indinf_{N_G(S)/S}^Gf_\un^{N_G(S)/S}\Defres_{N_G(S)/S}^G\,\varphi\mpoint
\end{equation}
For each $S\in\mathcal{G}$, the group $N_G(S)/S$ is a Roquette 2-group, i.e. it is cyclic, generalized quaternion of order at least 8, dihedral or semidihedral of order at least 16. Recall moreover that $\Im\,\imath_R=\Hom_\Z\big(R_\Q(R),\F_2\big)$ if $R$ is trivial, of order 2, or dihedral of order at least 16. If we show that $f_\un^{N_G(S)/S}\Defres_{N_G(S)/S}^G\,\varphi=0$ when $N_G(S)/S$ is not one of these Roquette 2-groups, then Equation~\ref{eq5} shows that $\varphi\in \Im\,\imath_G$, and we are done.\vspace{.5ex}\par
So all we have to do is to show that $f_\un^{N_G(S)/S}\Defres_{N_G(S)/S}^G\,\varphi=0$ when $R=N_G(S)/S$ is cyclic of order at least 4, generalized quaternion, or semi-dihedral. In each of these cases $R$ has a unique central subgroup $Z$ of order~2, and $f_\un^R=R/\un-R/Z\in B(R,R)$.\par
Let $X$ be any subgroup of $R$. Then
\begin{align*}
(f_\un^{N_G(S)/S}\Defres_{N_G(S)/S}^G\,\varphi)(R/X)&=(\Defres_{N_G(S)/S}^G\,\varphi)\big(f_\un^R\times_R(R/X)\big)\\
&=(\Defres_{N_G(S)/S}^G\,\varphi)(R/X-R/ZX)\mpoint
\end{align*}
This is zero if $Z\leq X$, so we can assume $Z\nleq X$, that is $Z\cap X=\un$.\par
If $R$ is cyclic or generalized quaternion, this implies $X=\un$. Then
\begin{align*}
(f_\un^{R}\Defres_{R}^G\,\varphi)(R/X)&=(\Defres_{R}^G\,\varphi)(R/\un-R/Z)\\
&=(\Defres_{R}^G\,\varphi)\big(\Ind_D^R(D/\un-D/Z)\big)\\
&=(\Res_D^R\Defres_{R}^G\varphi)(D/\un-D/Z)\\
&=(\Defres_D^G\varphi)(\varepsilon_D)\mvirg
\end{align*}
where $D$ is a cyclic subgroup of order 4 of $R$. By assumption, this is equal to $0_{\F_2}$, hence $f_\un^{N_G(S)/S}\Defres_{N_G(S)/S}^G\,\varphi=0$ in these cases, as was to be shown.\par
Now if $R\cong SD_{2^n}$ with $n\geq 4$, the same argument shows that 
$$(f_\un^{R}\Defres_{R}^G\,\varphi)(R/X)=0_{\F_2}\vspace{-1ex}$$
if $Z\leq X$ or $X=\un$. Up to conjugation, the only remaining subgroup $X$ of~$G$ is $X=I$, and then
\begin{align*}
(f_\un^{R}\Defres_{R}^G\,\varphi)(R/X)&=(\Defres_{R}^G\,\varphi)(R/I-R/IZ)\\
&=(\Defres_{R}^G\,\varphi)(\varepsilon_R)=0_{\F_2}\mpoint
\end{align*}
This completes the proof of Theorem~\ref{image}.\endpf
\begin{mth}{Corollary} Let $G$ be a finite group. Then the kernel of the natural map $\jmath_G:\F_2B(G)\to\Hom_{\F_2}\big(B^\times(G),\F_2\big)$ is the group generated by the elements $\Indinf_{T/S}^G\sur{\varepsilon}_{T/S}$, for all sections $(T,S)$ of $G$ such that $T/S\in\mathcal{R}$.
\end{mth}
\pf This is just a dual reformulation of Theorem~\ref{image}, since the natural map $\jmath_G:\F_2B(G)\to\Hom_{\F_2}\big(B^\times(G),\F_2\big)$ is the transposed map of $\imath_G$, up to the identification of $\F_2B(G)$ with the $\F_2$-bidual of $B(G)$.\endpf
For a finite group $R$, let $B_R=\Hom_{\CC}(R,-)=B(-,R)$ denote the representable biset functor defined by $R$, and let $\F_2B=\F_2\otimes_\Z B_R\cong \Hom_{\F_2\CC}(R,-)$.
\begin{mth}{Corollary} \label{reso}Let $\mathcal{S}$ be the set of finite groups defined by
$$\mathcal{S}=\{C_p\mid p\; \hbox{\rm odd prime}\}\cup\{C_4\}\cup\{SD_{2^n}\mid n\geq 4\}\mvirg\vspace{-1ex}$$
(by which we mean that $\CS$ contains exactly one group of order $p$ for each odd prime $p$, one cyclic group of order 4, and one semidihedral group of each order $2^n\geq 16$), 
and let $\CG_\mathcal{S}=\{(R,\sur{\varepsilon}_R)\mid R\in\mathcal{S}\}$ be the associated set of elements of $\F_2B$.
\begin{enumerate}
\item Let $L$ be the biset subfunctor of $\F_2B$ generated by $\CG_\CS$. Then there is an exact sequence of biset functors
$$\xymatrix{
0\ar[r]&L\ar[r]&\F_2B\ar[r]^-{\jmath}&\Hom_{\F_2}(B^\times,\F_2)\ar[r]&0\mpoint
}\vspace{-1ex}
$$
\item When $R\in\mathcal{R}$, let $d_R:\F_2B_R\to \F_2B$ be the morphism of biset functors induced by adjunction from $\sur{\varepsilon}_R\in\F_2B(R)$. Then the sequence
$$\xymatrix@C=3ex{
\mathop{\bigoplus}_{R\in\mathcal{S}}\limits\F_2B_R\ar[r]^-d&\F_2B\ar[r]^-\jmath&\Hom_{\F_2}(B^\times,\F_2)\ar[r]&0}\vspace{-1ex}
$$
is exact, where $d$ is the sum of all the maps $d_R$, for $R\in\mathcal{S}$.
\item For $R\in\CR$, let $\F_2B_Rf_\un^R$ be the direct summand of $\F_2B_R$ equal to the image of the idempotent endomorphism $f_\un^R\in B(R,R)=\End_{\CF}(B_R)$. Then the morphism $d_R:\F_2B_R\to \F_2B$ factors through a morphism $\hat{d}_R:\F_2B_Rf_\un^R\to\F_2B$, and this yields an exact sequence of functors
$$\xymatrix{
\mathop{\bigoplus}_{R\in \mathcal{S}}\limits\F_2B_Rf_\un^R\ar[r]^-{\hat{d}}&\F_2B\ar[r]^-\jmath&\Hom_{\F_2}(B^\times,\F_2)\ar[r]&0\mpoint}\vspace{-1ex}
$$
where $\hat{d}$ is the sum of the morphisms $\hat{d}_R$, for $R\in\mathcal{S}$.
\end{enumerate}
\end{mth}
\pf For Assertion 1, all we need to show is that the element $\sur{\varepsilon}_R$ belongs to $L(R)$, when $R$ is elementary abelian of order 4 or dihedral of order at least~8. But one checks easily that $\Res_{D_{2^{n-1}}}^{SD_{2^n}}\sur{\varepsilon}_{SD_{2^n}}=\sur{\varepsilon}_{D_{2^{n-1}}}$ for $n\geq 4$, and that $\Res_{(C_2)^2}^{D_8}\sur{\varepsilon}_{D_8}=\sur{\varepsilon}_{(C_2)^2}$. Now Assertion 2 follows from the fact that $L$ is equal to the image of $d$, and Assertion 3 from the fact that $f_\un^R\varepsilon_R=\varepsilon_R$ by Remark~\ref{invariant}.\endpf
\section{Extensions}\label{ext1}
\begin{mth}{Notation} \label{cloture}Let $\mathcal{T}$ denote the class of finite groups which are subquotient of some group in $\mathcal{S}$, that is the class of groups isomorphic to a group in the following set
$$\left\{\begin{array}{c}\{C_p\mid p\;\hbox{\rm odd prime}\}\cup\{C_{2^n}\mid n\geq 0\}\cup\{Q_{2^n}\mid n\geq3\}\\\rule{0ex}{3ex}\cup\{(C_2)^2\}\cup\{D_{2^n}\mid n\geq 3\}\cup\{SD_{2^n}\mid n\geq 4\}\end{array}\right.$$
\end{mth}
In other words $\CT$ is the class of groups which are cyclic of odd prime order, elementary abelian of order 4, dihedral of order 8, or a Roquette 2-group. 
\begin{mth}{Proposition} \begin{enumerate}
\item Set $\widehat{B^\times}=\Hom_{\F_2}(B^\times,\F_2)$. If $S$ be a biset functor over $\F_2$, such that $S(\un)=\zero$, then $\Ext^1_{\CF_{\F_2}}(\widehat{B^\times},S)$ is isomorphic to a subspace of $\mathop{\prod}_{R\in\mathcal{S}}\limits\partial S(R)$, where $\mathcal{S}$ is the set defined in Corollary~\ref{reso}.
\item Let $H$ be a finite group. If there exists a simple $\F_2\Out(H)$-module $V$ such that $\Ext^1_{\CF_{\F_2}}(S_{H,V},B^\times)\neq \zero$, then $H\in\CT$. 
\end{enumerate}
\end{mth}
\pf The functor $L$ of Corollary~\ref{reso} is equal to $\Im\,d=\Im\hat{d}=\Ker\,\jmath$, so we have an exact sequence
$$\xymatrix{0\ar[r]&L\ar[r]&\F_2B\ar[r]^-\jmath&\widehat{B^\times}\ar[r]&0\mvirg}
$$
and a surjective morphism $\xymatrix{\mathop{\bigoplus}_{R\in \mathcal{S}}\limits\F_2B_Rf_\un^R\ar@{->>}[r]^-{\hat{d}}&L}$. Let $M$ be any biset functor over $\F_2$. Applying first the functor $\Hom_{\CF_{\F_2}}(-,M)$ to the short exact sequence, we get the beginning of a long exact sequence
$$\xymatrix@C=2.5ex{0\ar[r]&\Hom(\widehat{B^\times},M)\ar[r]&\Hom(\F_2B,M)\ar[r]&\Hom(L,M)\ar[r]&\Ext^1(\widehat{B^\times},M)\ar[r]&0\mvirg}$$
where $\Hom$ and $\Ext^1$ are taken in the category $\CF_{\F_2}$. Indeed $\Ext^1(\F_2B,M)$ is equal to zero as $\F_2B=\Hom_{\CF_{\F_2}}(\un,-)$ is representable, hence projective in $\CF_{\F_2}$. Moreover $\Hom(\F_2B,M)\cong M(\un)$, so we get an exact sequence
\begin{equation}\label{ext}\xymatrix@C=2.5ex{0\ar[r]&\Hom(\widehat{B^\times},M)\ar[r]&M(\un)\ar[r]&\Hom(L,M)\ar[r]&\Ext^1(\widehat{B^\times},M)\ar[r]&0\mpoint}
\end{equation}
In the case $M=S$, since $S(\un)=\zero$, this gives an isomorphism $\Hom(L,S)\cong\Ext^1(\widehat{B^\times},S)$. On the other hand, since $L$ is a quotient of$\mathop{\bigoplus}_{R\in \mathcal{S}}\limits\F_2B_Rf_\un^R$, we get an inclusion
$$\xymatrix{\Hom(L,S)\ar@{^{(}->}[r]&\Hom(\mathop{\bigoplus}_{R\in \mathcal{S}}\limits\F_2B_Rf_\un^R,S)\cong\prod_{R\in\mathcal{S}}\limits\Hom(\F_2B_Rf_\un^R,S)\mpoint}$$
Moreover $\Hom(\F_2B_Rf_\un^R,S)\cong f_\un^RS(R)=\partial S(R)$. This proves Assertion 1.\par
For Assertion 2, we take $M=S_{H,W}$ in the exact sequence~\ref{ext}, where $W=V^*$ is the dual module. If $H\neq\un$, then $S_{H,W}(\un)=\zero$. And if $H=\un$, then $W=\F_2$, and by duality
$$\Hom(\widehat{B^\times},S_{\un,\F_2})\cong \Hom(S_{\un,\F_2},B^\times)\cong\F_2\cong S_{\un,\F_2}(\un)\mpoint$$
In both cases, we get an isomorphism  
$$\Hom(L,S_{H,W})\cong\Ext^1(\widehat{B^\times},S_{H,W})\mpoint$$ 
Now by duality, this is isomorphic to $\Ext^1(S_{H,V},B^\times)$. As before, it embeds into $\prod_{R\in\mathcal{S}}\limits \partial S_{H,W}(R)$. Hence if $\Ext^1(S_{H,V},B^\times)\neq\zero$, then there exists some $R$ in $\mathcal{S}$ such that $\partial S_{H,W}(R)\neq\zero$. In particular $S_{H,W}(R)\neq\zero$, so $H$ is a subquotient of $R$. This completes the proof.\endpf
\section{Sectional characterization}\label{char}
Let $\CO_\CT$ denote the forgetful functor $\CF_k\to\CF_{\CT,k}$. It was shown in Section~5 of~\cite{both3} that the functor $\CO_\CT$ has a right adjoint $\CR_\CT$ defined as follows\footnote{The construction in~\cite{both3} actually dealt only with $p$-biset functors, for a fixed prime number $p$, and their restriction to a class of finite $p$-groups closed under taking subquotients. But it extends {\em verbatim} to the categories $\CF_k$ and $\CF_{\CT,k}$, for any class $\CT$ of finite groups closed under taking subquotients.}. For a finite group $G$, let $\CT(G)$ denote the set of sections $(T,S)$ of $G$ such that $T/S\in\CT$. If $F\in\CF_{\CT,k}$, then
$$\CR_\CT(F)(G)=\limproj{(T,S)\in\CT(G)}F(T/S)\mvirg$$
that is the set of sequences of elements $l_{T,S}\in F(T/S)$, for $(T,S)\in\CT(G)$, subject to the following conditions:
\begin{enumerate}
\item If $(T,S)$ and $(T',S')$ are elements of $\CT(G)$ such that $S\leq S'\leq T'\leq T$, then
$$\Defres_{T'/S'}^{T/S}l_{T,S}=l_{T',S'}\mvirg$$
\item If $(T,S)\in\CT(G)$ and $x\in G$, then
$$^xl_{T,S}=l_{^xT,{^xS}}\mpoint$$
\end{enumerate}
The biset functor structure on $\CR_\CT(F)$ is obtained as follows. For a finite group $H$ and a finite $(H,G)$-biset $U$, the image $v=U\,l$ of $l\in\CR_\CT(F)(G)$ by~$U$ is the sequence $m_{T,S}$, for $(T,S)\in\CT(H)$, defined by
$$m_{T,S}=\sum_{u\in[T\dom U/G]}(S\dom Tu)\,l_{T^u,S^u}\mvirg$$
where for a subgroup $X$ of $H$, we set $X^u=\{g\in G\mid\exists x\in X,\;xu=ug\}$. This makes sense as one can show that $(T^u,S^u)$ is a section of $G$, and that $T^u/S^u$ is isomorphic to a subquotient of $T/S$, hence in $\CT$ if $T/S\in \CT$.\par
The unit $\eta:\Id\to \CR_\CT\circ\CO_\CT$ of the adjunction between $\CO_\CT$ and $\CR_\CT$, evaluated for $N\in \CF_k$ at a finite group $G$, is the map $\eta_{N,G}$ from $N(G)$ to $\CR_\CT\CO_\CT(N)(G)$ sending $n\in N(G)$ to the sequence of elements $l_{T,S}\in \CO_\CT(N)(T/S)=N(T/S)$, for $(T,S)\in \CT(G)$, defined by
$$l_{T,S}=\Defres_{T/S}^Gn\mpoint$$
\begin{mth}{Theorem} The unit of the adjunction between $\CO_\CT:\CF_{\F_2}\to \CF_{\CT,\F_2}$ and $\CR_\CT:\CF_{\CT,\F_2}\to\CF_{\F_2}$ induces an isomorphism of biset functors over $\F_2$
$$\eta_{B^\times}:B^\times\to \CR_\CT\CO_\CT(B^\times)\mpoint$$
In other words, for any finite group $G$, the map
$$\eta_{B^\times,G}:B^\times(G)\to \limproj{(T,S)\in\CT(G)}B^\times(T/S)$$
is an isomorphism.
\end{mth}
\pf Set $\widehat{B}=\Hom_\Z(B,\F_2)$, and let $C$ denote the cokernel of the morphism $\imath: B^\times\to \widehat{B}$. The functor $\CO_T$ is exact, so its right adjoint $\CR_T$ is left exact. Hence we get a commutative diagram
\begin{equation}\label{serpent}\xymatrix{0\ar[r]&B^\times\ar[r]^-\imath\ar[d]^-{\eta_{B^\times}}&\widehat{B}\ar[r]^-\pi\ar[d]^-{\eta_{\widehat{B}}}&C\ar[r]\ar[d]^-{\eta_{C}}&0\\
            0\ar[r]&\CR_\CT\CO_\CT(B^\times)\ar[r]^-f&\CR_\CT\CO_\CT(\widehat{B})\ar[r]^-g&\CR_\CT\CO_\CT(C)
}
\end{equation}
with exact rows, where $\pi$ is the projection morphism, $f=\CR_\CT\CO_\CT(\imath)$, and $g=\CR_\CT\CO_\CT(\pi)$.\par
We first claim that the vertical morphism $\eta_{\widehat{B}}$ is an isomorphism. Indeed, let $G$ be a finite group. The map $\eta_{\widehat{B},G}:\widehat{B}(G)\to\limproj{(T,S)\in \CT(G)}\widehat{B}(T/S)$ is the map sending the linear form $\varphi:B(G)\to\F_2$ to the sequence of linear forms $d_{T,S}=\Defres_{T/S}^G\varphi:B(T/S)\to\F_2$, for $(T,S)\in\CT(G)$. Since
$$d_{T,S}\big((T/S)\big/(T/S)\big)=\varphi\big(\Indinf_{T/S}^G(T/S)\big/(T/S)\big)=\varphi(G/T)\mvirg$$
and since $(T,T)\in\CT$ for any subgroup $T$ of $G$, we see that $\eta_{\widehat{B},G}$ is injective.\par
Now let $l=(\psi_{T,S})_{(T,S)\in\CT(G)}$ be an element of $\limproj{(T,S)\in\CT(G)} \widehat{B}(T/S)$. So for each $(T,S)\in\CT(G)$, we get a linear map $\psi_{T,S}:B(T/S)\to\F_2$. For a subgroup $U$ of $G$, we set $\varphi(G/U)=\psi_{U,U}\big((U/U)\big/(U/U)\big)$. Then for any $x\in G$, we have $\varphi(G/{^xU})=\varphi(G/U)$, since
\begin{align*}
\varphi(G/{^xU})&=\psi_{^xU,^xU}\big((^xU/^xU)\big/(^xU/^xU)\big)\\
&=(^x\psi_{U,U})\big((^xU/^xU)\big/(^xU/^xU)\big)\\
&=\psi_{U,U}\big((U/U)\big/(U/U)\big)
\end{align*}
because $^x\psi_{T,S}=\psi_{^xT,^xS}$ for any $(T,S)\in\CT(G)$. It follows that $\varphi$ extends linearly to an element of $\widehat{B}(G)$. Moreover, for any $(T,S)\in\CT(G)$ and any subgroup $U/S$ of $T/S$, we have $\psi_{U,U}=\Defres_{U/U}^{T/S}\psi_{T,S}$, so
\begin{align*}
\Defres_{T/S}^G\varphi\big((T/S)\big/(U/S)\big)&=\varphi\Big( \Indinf_{T/S}^G\big((T/S)\big/(U/S)\big)\Big)\\
&=\varphi(G/U)=\psi_{U,U}\big((U/U)\big/(U/U)\big)\\
&=(\Defres_{U/U}^{T/S}\psi_{T,S})\big((U/U)\big/(U/U)\big)\\
&=\psi_{T,S}\Big(\Indinf_{U/U}^{T/S}\big((U/U)\big/(U/U)\big)\\
&=\psi_{T,S}\big((T/S)\big/(U/S)\big)\mpoint
\end{align*}
Hence $\Defres_{T/S}^G\varphi=\psi_{T,S}$, for any $(T,S)\in\CT(G)$, so the map $\eta_{\widehat{B},G}$ is surjective. Hence it is an isomorphism, for any finite group $G$, and $\eta_{\widehat{B}}$ is an isomorphism, as claimed.\par
Now the Snake's lemma, applied to Diagram~\ref{serpent}, shows that the morphism $\eta_{B^\times}$ is injective, and that its cokernel is isomorphic to the kernel of $\eta_C$. For a finite group $G$, an element $\sur{\varphi}$ of $\Ker\,\eta_{C,G}$ is represented by a linear form $\varphi:B(G)\to\F_2$ such that $\Defres_{T/S}^G\varphi$ belongs to $\imath\big(B^\times(T/S)\big)$, for any $(T,S)\in \CT(G)$, hence in particular for any section $(T,S)$ of $G$ such that $T/S\in\mathcal{R}$. By Theorem~\ref{image}, it follows that $(\Defres_{T/S}^G\varphi)(\epsilon_{T,S})=0$. As this holds for any section $(T,S)$ of $G$ with $T/S\in\CR$, the form $\varphi$ lies in the image of $\imath_G$, by Theorem~\ref{image} again. In other words the element $\sur{\varphi}$ is equal to zero. Hence $\Ker\,\eta_{C,G}=\zero$ for any $G$, so $\eta_C$ is injective. It follows that $\eta_{B^\times}$ is surjective. Hence it is an isomorphism. This completes the proof of the theorem.\endpf
\section{Finite generation and presentation}\label{pres}
In this section, we show that the biset functor $B^\times$ is not finitely generated, and that its dual $\widehat{B^\times}$ is finitely generated (by a single element!), but not finitely presented. Recall that $k$ is a commutative ring, and that $kB_R=k\otimes_\Z B(-,R)$ is the representable functor $\Hom_{k\CC}(R,-)$.
\begin{mth}{Lemma} \label{residue}Let $R$ and $H$ be finite groups. If $\sur{kB_R}(H)\neq\zero$, then $H$ is a subquotient of $R$.
\end{mth}
\pf Recall that $kB_R(H)=kB(H,R)$ has a $k$-basis consisting of the transitive bisets $(H\times R)/X$, for $X$ in a set of representatives of conjugacy classes of subgroups of $H\times R$. Let $X$ be one of these groups, and $Y=p_1(X)$ be its first projection. Then $(H\times R)/X=kB_R(\Ind_Y^H)\big((Y\times R)/X\big)$, where we abuse notation in the right hand side by viewing $X$ as a subgroup of $Y\times R$, and $(Y\times R)/X$ as an element of $kB_R(Y)$. Hence if the image of $(H\times R)/X$ in $\sur{kB_R}(H)$ is non zero, then $Y=H$. \par
In this case let $N=k_1(X)$. Then $N$ is a normal subgroup of $H=p_1(X)$, and $(H\times R)/X=kB_R(\Inf_{H/N}^H)\Big(\big((H/N)\times R\big)/\sur{X}\Big)$, where 
$$\sur{X}=\{(hN,g)\mid (h,g)\in X\}\mpoint$$
Hence if the image of $(H\times R)/X$ in $\sur{kB_R}(H)$ is non zero, we also have $N=\un$. Then $H\cong Y/N\cong T/S$, where $T=p_2(X)$ and $S=k_2(X)$. In particular $H$ is a subquotient of $R$. This completes the proof.\endpf
\begin{mth}{Proposition} \label{fg}Let $k$ be a commutative ring and $\CT$ be a class of finite groups closed under taking subquotients. Let $F\in\CF_{\CT,k}$. The following conditions are equivalent:
\begin{enumerate}
\item The functor $F$ is finitely generated.
\item There exists a finite family $\mathcal{E}$ of groups in $\CT$ and an epimorphism $\mathop{\oplus}_{R\in\mathcal{E}}\limits kB_R\to F$.
\item For any $H\in \CT$, the $k$-module $\sur{F}(H)$ is finitely generated, and there exists an integer $n\in \N$ such that $\sur{F}(H)=\zero$ whenever $|H|>n$.
\end{enumerate}
\end{mth}
\pf The equivalence of 1 and 2 is classical. If 1 holds, then there is a finite set $\mathcal{G}$ of elements of $F$ such that $\langle\CG\rangle=F$. For each $(R,v)\in\CG$, we get a morphism $\widetilde{v}:kB_R\to F$ associated to $v\in F(R)$ by Yoneda's lemma. The sum of these morphisms
$$\mathop{\oplus}_{(R,v)\in\CG}\widetilde{v}:\mathop{\oplus}_{(R,v)\in\CG}\limits kB_R\to F$$
is surjective if $\langle\CG\rangle=F$, so 2 holds. Conversely, each representable functor $kB_R$ is generated by the single element $(R,\Id_R)$, where $\Id_R\in kB(R,R)$ is the identity endomorphism of $R$. Hence if 2 holds, then $F$ is a quotient of a finite sum of finitely generated functors, so $F$ is finitely generated, and 1 holds.\medskip\par
Now if 2 holds, then for each $H\in\CT$, the $k$-module $F(H)$ is a quotient of $\mathop{\oplus}_{R\in\mathcal{E}}\limits kB(H,R)$, and each $kB(H,R)$ is a finitely generated $k$-module. Hence $F(H)$ is a finitely generated $k$-module. Moreover $\sur{F}(H)$ is a quotient of $\mathop{\oplus}_{R\in\mathcal{E}}\limits \sur{kB_R}(H)$, and $\sur{kB_R}(H)=\zero$ unless $H$ is a subquotient of~$R$, by Lemma~\ref{residue}. In particular $\sur{kB_R}(H)=\zero$ if $|H|>|R|$, so $\sur{F}(H)=\zero$ if $|H|>n$, where $n=\max\{|R|\mid R\in \CE\}$. Hence 3 holds.\medskip\par
Now if 3 holds, there is a finite set of isomorphism classes of finite groups $R$ such that $\sur{F}(R)\neq \zero$. Let $\mathcal{U}$ be a set of representatives of this set. For each $R\in\mathcal{U}$, we can lift to $F(R)$ a finite generating set of the $k$-module $\sur{F}(R)$. We get a finite subset $V_R$ of $F(R)$, and this gives a finite set 
$$\CG=\{(R,v)\mid R\in\mathcal{U},\;v\in V_R\}$$
of elements of $F$, which in turns gives a morphism
$$\pi:P=\mathop{\oplus}_{(R,v)\in\CG}kB_R\to F\mpoint$$
Our choice of $\mathcal{U}$ and $V_R$, for $R\in\mathcal{U}$, shows that the induced morphism $\sur{\pi_H}:\sur{P}(H)\to \sur{F}(H)$ is surjective for any $H\in\CT$: if $H$ is not isomorphic to a group in $\mathcal{U}$, then this is trivially true because $\sur{F}(H)=\zero$. And otherwise, we can assume $H\in\mathcal{U}$, and then $\sur{\pi_H}$ is surjective because $\sur{\pi_H}(V_H)$ generates $\sur{F}(H)$ by construction.\par
We deduce by induction on $n=|H|$ that $\pi_H:P(H)\to F(H)$ is surjective for any $H\in\CT$. For $n=1$, this is clear, since 
$$\pi_\un=\sur{\pi_1}:\sur{P}(\un)=P(\un)\to F(\un)=\sur{F}(\un)$$
is surjective. Now assume $\pi_K$ is surjective for any $K\in \CT$ with $|K|<n=|H|$, and let $v\in F(H)$. Since $\sur{\pi_H}$ is surjective, there is an element $w\in P(H)$, a set $\Sigma$ of proper sections of $H$ (i.e. sections different from $(H,\un)$), and elements $v_{T,S}\in F(T/S)$, for $(T,S)\in\Sigma$, such that
$$v=\pi_H(w)+\sum_{(T,S)\in\Sigma}\Indinf_{T/S}^Hv_{T,S}\mpoint$$
Since $|T/S|<n$ for any $(T,S)\in \Sigma$, the map $\pi_{T/S}:P(T/S)\to F(T/S)$ is surjective, and there is an element $w_{T,S}\in P(T/S)$ such that $v_{T/S}=\pi_{T/S}(w_{T,S})$. It follows that
\begin{align*}
v&=\pi_H(w)+\sum_{(T,S)\in\Sigma}\Indinf_{T/S}^H\pi_{T/S}(w_{T,S})\\
&=\pi_H(w)+\sum_{(T,S)\in\Sigma}\pi_H\big(\Indinf_{T/S}^Hw_{T,S}\big)=\pi_H\big(w+\sum_{(T,S)\in\Sigma}\Indinf_{T/S}^Hw_{T,S}\big)\mpoint\\
\end{align*}
Hence $\pi_H$ is surjective, and this completes the inductive step.\par
It follows that $\pi$ is an epimorphism, so 3 implies 2, completing the proof of Proposition~\ref{fg}.\endpf
\begin{mth}{Corollary} The biset functor $B^\times$ is not finitely generated.
\end{mth}
\pf It has been show by Barsotti (\cite{barsotti}, Proposition 6.8) that if $p$ is a prime number congruent to 1 mod 4, then $\sur{B^\times}(D_{2p})\neq \zero$, where $D_{2p}$ is a dihedral group of order $2p$ (in Barsotti's terminology, the group $D_{2p}$ is {\em residual}). It follows that there are arbitrary large finite groups $H$ such that  $\sur{B^\times}(H)\neq\zero$. By Proposition~\ref{fg}, the functor $B^\times$ is not finitely generated.\endpf
Recall from Corollary~\ref{reso} that there is an exact sequence of biset functors
\begin{equation}\label{ses}\xymatrix{
0\ar[r]&L\ar[r]&\F_2B\ar[r]^-{\jmath}&\widehat{B^\times}\ar[r]&0\mvirg
}
\end{equation}
where 
$$\mathcal{S}=\{C_p\mid p\; \hbox{\rm odd prime}\}\cup\{C_4\}\cup\{SD_{2^n}\mid n\geq 4\}\mvirg$$
and $L$ is the biset subfunctor of $\F_2B$ generated $\CG_\CS=\{(R,\sur{\varepsilon}_R)\mid R\in\CS\}$. 

\begin{mth}{Proposition} \begin{enumerate}
\item Let $\CG=\{(\un,u)\}$, where $u$ is the non zero element of $\widehat{B^\times}(\un)\cong\F_2$. Then $\langle\CG\rangle=\widehat{B^\times}$.
\item The functor $L$ is not finitely generated.
\item The functor $\widehat{B^\times}$ is not finitely presented.
\end{enumerate}
\end{mth}
\pf (1) This follows from the fact that $\widehat{B^\times}$ is a quotient of $\F_2B$, and that $\F_2B$ is generated by ${(\un,e)}$, where $e\in\F_2B(\un)$ is the class of a set of cardinality one, endowed with the trivial action of the trivial group. Indeed if $K\leq H$ are finite groups, then $H/K=\Indinf^H_{K/K}\Iso(f_K)(e)$, where $f_K:\un\to K/K$ is the unique group isomorphism.\mpn
(2) The exact sequence~(\ref{ses}) shows that if $p$ is an odd prime number, then $L(C_p)\cong \F_2$ and $L(\un)=\zero$. Hence $\sur{L}(C_p)\cong L(C_p)\cong \F_2$, so there exist arbitrary large finite groups $H$ such that $\sur{L}(H)\neq\zero$. By Proposition~\ref{fg}, the functor $L$ is not finitely generated.\mpn
(3) Suppose that there exists an exact sequence in $\CF_{\F_2}$
\begin{equation}\label{ses2}\xymatrix{N\ar[r]&M\ar[r]&\widehat{B^\times}\ar[r]&0\mvirg}
\end{equation}
where $M$ is projective and $N$ is finitely generated. This gives an exact sequence
$$\xymatrix{0\ar[r]&K\ar[r]&M\ar[r]&\widehat{B^\times}\ar[r]&0\mvirg}$$
where $K$ is the image of $N$ in $M$. In particular $K$ is finitely generated. Then, since $M$ and $\F_2B$ are projective in $\CF_{\F_2}$, Shanuel's lemma gives an isomorphism of functors
$$L\oplus M\cong K\oplus \F_2B\mpoint$$
Then $L$ is a quotient of $K\oplus \F_2B$, which is finitely generated. Hence $L$ is finitely generated, contradicting 2. So no short exact sequence like~(\ref{ses2}) can exists, hence $\widehat{B^\times}$ is not finitely presented.\endpf
Recall that in Corollary~\ref{reso}, the set
$$\mathcal{S}=\{C_p\mid p\; \hbox{\rm odd prime}\}\cup\{C_4\}\cup\{SD_{2^n}\mid n\geq 4\}$$
was introduced. We finally prove that $\CG_\CS=\{(R,\sur{\varepsilon}_R)\mid R\in \CS\}$ is a minimal set of generators of $L$.
\begin{mth}{Theorem} \label{minimal}Let $\CS'$ be a proper subset of the set $\CS$ introduced in Corollary~\ref{reso}, and $\CG_{\CS'}=\{(R,\sur{\varepsilon}_R)\mid r\in\CS'\}$. Then $\langle\CG_{\CS'}\rangle$ is a proper subfunctor of $\langle\CG_{\CS}\rangle=L$.
\end{mth}
\pf It suffices to show that for any $R\in\CS$, if we set $\CG_\CS^R=\CG_\CS-\{(R,\sur{\varepsilon}_R)\}$ and $L^R=\langle \CG_\CS^R\rangle$, then $\sur{\varepsilon}_R\notin L^R(R)$. So we assume that $\sur{\varepsilon}_R\in L^R(R)$, for contradiction, i.e. that there exists a finite set $\CE$ of pairs $(H,X)$, where $H\in \CS-\{R\}$ and $X$ is a subgroup of $R\times H$, such that
\begin{equation}\label{epsilon}
\sur{\varepsilon}_R=\sum_{(H,X)\in\CE}\big((R\times H)/X\big)\sur{\varepsilon}_H\mpoint
\end{equation}
Since $f_\un^R\sur{\varepsilon}_R=\sur{\varepsilon}_R$ for any $R$, by Remark~\ref{invariant}, this also reads
$$\sur{\varepsilon}_R=\sum_{(H,X)\in\CE}f_\un^R\big((R\times H)/X\big)f_\un^H\sur{\varepsilon}_H\mpoint$$
We now observe that each element $H$ of $\CS$ has a unique minimal normal (hence central of prime order) subgroup $Z_H$. And if $k_1(X)\geq Z_R$, then 
$$f_\un^R\big((R\times H)/X\big)=f_\un^R\Inf_{R/Z_R}^R\Def_{R/Z_R}^R\big((R\times H)/X\big)=0$$
by Lemma 6.3.2 of \cite{bisetfunctorsMSC}. Similarly $\big((R\times H)/X\big)f_\un^H=0$ if $k_2(X)\geq Z_H$. So we can assume that $k_1(X)\cap Z_R=\un$ and $k_2(X)\cap Z_H=\un$ for any $(H,X)\in\CE$. \par
This forces $k_1(X)=k_2(X)=\un$, unless $R$ is semidihedral and $k_1(X)$ is a non central subgroup of order 2 of $R$, or $H$ is semidihedral and $k_2(X)$ is a non central subgroup of order 2 of $H$ (both cases may occur simultaneously). In the first case $p_1(X)\leq N_R\big(k_1(X)\big)=k_1(X)Z_R$, and $p_1(X)/k_1(X)$ has order 1 or 2. Similarly, in the second case $p_2(X)/k_2(X)$ has order 1 or 2. In any of these two cases, the morphism $f_\un^R\big((R\times H)/X\big)f_\un^H$ of $\CC_{\F_2}$ factors through a group of order 1 or 2 (see~\ref{decomposition}). Since $L(\un)=L(C_2)=\zero$, it follows that $f_\un^R\big((R\times H)/X\big)f_\un^H\sur{\varepsilon}_H=0$. \par
So we can assume that $k_1(X)=k_2(X)=\un$ for any $(H,X)\in\CE$. In this case $X$ is a twisted diagonal subgroup of $R\times H$, that is, there is a subgroup $K$ of $H$, a subgroup $S$ of $R$, and a group isomorphism $f:K\to S$, such that 
\begin{equation}\label{induit}
(R\times H)/X\cong \Ind_S^R\,\Iso(f)\,\Res_K^H\mpoint
\end{equation}
We can also assume that $\Res_K^H\sur{\varepsilon}_H\neq 0$, and in particular that $L(K)\neq\zero$.\smallskip\par
Suppose first that $R=C_p$ for an odd prime number $p$. Then $S=\un$ or $S=R$. Since $L(\un)=\zero$, we have $S=R$. Then $K\cong S=C_p$ is a subgroup of $H\in \CS$. The only element $H$ of $\CS$ admitting a subgroup of odd prime order $p$ is $C_p$ itself. But $H\in\CS-\{C_p\}$ by assumption, so we get a contradiction. \smallskip\par
Suppose now that $R=C_4$. Again, since $L(\un)=L(C_2)=\zero$, we have $S=R$, and $H\geq K\cong C_4$. Since $H\in \CS-\{C_4\}$, it follows that $H=SD_{2^n}$ for some $n\geq 4$. Then $K$ is contained in the unique generalized quaternion subgroup $Q$ of index 2 of $H$. One checks easily that $\Res_Q^H\sur{\varepsilon}_H=Q/\un-Q/Z_H$, and it follows that $\Res_K^H\sur{\varepsilon}_H=0$. We get a contradiction also in this case.\smallskip\par
We are left with the case $R=SD_{2^n}$, for some $n\geq 4$. Then if $(H,X)\in\CE$, the first projection $S$ of the twisted diagonal subgroup $X$ of $R\times H$ is a proper subgroup of $R$: otherwise indeed, its second projection $K$ is a semidihedral subgroup of $H\in\CS$, which can only occur if $H$ itself is semidihedral and isomorphic, hence equal, to $R$. This is a contradiction, since $H\in\CS-\{R\}$.\par
It follows from (\ref{epsilon}) and (\ref{induit}) that $\sur{\varepsilon}_R$ is a sum of elements of the form $\Ind_S^Ru_S$, for proper subgroups $S$ of $R$ and elements $u_S$ of $L(S)$. Let $D$, $C$, and $Q$ be the subgroups of index 2 of $R$, where $D$ is dihedral, $C$ is cyclic, and $Q$ is generalized quaternion. We can write
$$\sur{\varepsilon}_R=\Ind_D^Rv_D+\Ind_C^Rv_C+\Ind_Q^Rv_Q\mvirg$$
for some $v_D\in L(D)$, $v_C\in L(C)$, and $v_Q\in L(Q)$. We set $Z=Z_R$ for simplicity. If $M\in\{D,C,Q\}$, then $M\geq Z$, and one checks easily that $f_\un^R\Ind_M^R=\Ind_M^Rf_\un^M$. It follows that
$$\sur{\varepsilon}_R=\Ind_D^Rf_\un^Dv_D+\Ind_C^Rf_\un^Cv_C+\Ind_Q^Rf_\un^Qv_Q\mpoint$$
Moreover $f_\un^MM/N=0$ if $N\cap Z(M)\neq \un$. It follows that $f_\un^Cv_C$ is a multiple of $C/\un-C/Z$, and that $f_\un^Qv_Q$ is a multiple of $Q/\un-Q/Z$. Hence
$$\sur{\varepsilon}_R=\Ind_D^Rw_D+\lambda(R/\un-R/Z)\mvirg$$
for some $\lambda\in \F_2$ and some $w_D(=f_\un^Dv_D)\in \partial L(D)$.\par
Now cutting the exact sequence
$$\xymatrix{
0\ar[r]&L(D)\ar[r]&\F_2B(D)\ar[r]&\widehat{B^\times}(D)\ar[r]&0}
$$
by the idempotent $f_\un^D$ gives the exact sequence
$$\xymatrix{
0\ar[r]&\partial L(D)\ar[r]&\partial \F_2B(D)\ar[r]&\partial\widehat{B^\times}(D)\ar[r]&0\mpoint}
$$
The vector space $\partial \F_2B(D)$ has a basis consisting of the elements $D/\un-D/Z$, $D/I-D/IZ$, and $D/J-D/JZ$, where $I$ and $J$ are two non conjugate non central subgroups of order 2 of $D$. The vector space $\partial\widehat{B^\times}(D)$ is isomorphic to the dual of $\partial B^\times(D)$, which is one dimensional (see e.g. Corollary~6.12 of~\cite{burnsideunits}, or the proof of Theorem~\ref{image}). It follows that $\partial L(D)$ has dimension 2. Now $\partial L(D)$ contains the two elements $\sur{\varepsilon}_D=(D/I-D/IZ)-(D/J-D/JZ)$ and $D/\un-D/Z=\Ind_{C_4}^D\sur{\varepsilon}_{C_4}$, which are obviously linearly independent. Hence these two elements form a basis of $\partial L(D)$. It follows that
$$\sur{\varepsilon}_R=\alpha\Ind_D^R\sur{\varepsilon}_D+\beta\Ind_{C_4}^R\sur{\varepsilon}_{C_4}+\lambda(R/\un-R/Z)\mvirg$$
for some $\alpha, \beta$ in $\F_2$. But one checks easily that $\Ind_D^R\sur{\varepsilon}_D=0$, for $I$ and $J$ are conjugate in $R$. Moreover $\Ind_{C_4}^R\sur{\varepsilon}_{C_4}=R/\un-R/Z$. Thus $\sur{\varepsilon}_R=R/I-R/IZ$ is a scalar multiple of $R/\un-R/Z$, which is obviously wrong. This final contradiction completes the proof of Theorem~\ref{minimal}.\endpf
\begin{mth}{Corollary} Let $2^\CS$ be the set of subsets of $\CS$, ordered by inclusion of subsets, and $[0,L]$ be the poset of subfunctors of $L$, ordered by inclusion of subfunctors. Let $g:2^\CS\to [0,L]$ be the map sending $A\subseteq\CS$ to 
$$g(A)=\big\langle \{(H,\sur{\varepsilon}_H)\mid H\in A\}\big\rangle\subseteq L\mvirg$$
and $f:[0,L]\to 2^\CS$ be the map sending the subfunctor $M$ of $L$ to
$$f(M)=\{H\in\CS\mid \sur{\varepsilon}_H\in M(H)\}\subseteq \CS\mpoint$$
Then:
\begin{enumerate}
\item Let $A, A'$ be subsets of $\CS$, and $M, M'$ be subfunctors of $L$. Then 
$$g(A\cup A')=g(A)+g(A')\;\hbox{and}\;f(M\cap M')=f(M)\cap f(M')\mpoint$$
In particular $f$ and $g$ are maps of posets.
\item $f\circ g=\Id_{2^\CS}$. 
\item The poset $[0,L]$ is uncountable.
\end{enumerate}
\end{mth}
\pf Assertions 1 is straightforward. For Assertion 2, let $A\subseteq \CS$ and $A'=f\circ g(A)$. Then clearly $A\subseteq A'$. If this inclusion is strict, let $S\in A'-A$. Then 
$$(S,\sur{\varepsilon}_S)\in g(A)=\big\langle \{(H,\sur{\varepsilon}_H)\mid H\in A\}\big\rangle\subseteq g(\CS')\mvirg$$
where $\CS'=\CS-\{S\}$. Then $g(\CS')=g(\CS)=L$, and by Theorem~\ref{minimal}, it follows that $\CS'=\CS$, a contradiction. Hence $A=A'$, and $f\circ g=\Id_{2^\CS}$. In particular $g$ is injective, and Assertion 3 follows, since the set of subsets of the (infinite) countable set $\CS$ is uncountable.\findemo
\begin{rem}{Remark} The map $g$ is not surjective, and not a map of lattices (that is, the image by $g$ of an intersection of subsets need not be the intersection of the images of the subsets). Indeed, if $S$ is any semidihedral group in $\CS$, and $C$ its cyclic subgroup of index 2, we have
$$\Res_C^S\sur{\varepsilon}_S=C/\un-C/Z=\Ind_{C'}^C\sur{\varepsilon}_{C'}\mvirg$$
where $Z$ is the center of $S$ and $C'$ the subgroup of order 4 of $C$. It follows that $u=C/\un-C/Z$ is a non-zero element of $M(C)$, where $M$ is the intersection of the subfunctors of $L$ generated by $\{(S,\sur{\varepsilon}_S)\}$ and $\{(C_4,\sur{\varepsilon}_{C_4})\}$. 
In other words if $A=\{S\}$ and $A'=\{C_4\}$, we have $0\neq u\in \big(g(A)\cap g(A')\big)(C)$, so $g(A)\cap g(A')\neq\zero$. Now if $g(A)\cap g(A')$ belongs to the image of $g$, there is a subset $A''$ of $\CS$ such that $g(A)\cap g(A')=g(A'')$, and then
$$A''=fg(A'')=f\big(g(A)\cap g(A')\big)=fg(A)\cap fg(A')=A\cap A'=\emptyset\mpoint$$
This is a contradiction since $g(A'')\neq\zero$ but $g(\emptyset)=\zero$. It follows that $g$ is not surjective, and not a map of lattices, since $g(A)\cap g(A')\neq g(A\cap A')$. 
\end{rem}

\begin{thebibliography}{10}

\bibitem{barsotti}
J.~Barsotti.
\newblock On the unit group of the {B}urnside ring as a biset functor for some
  solvable groups.
\newblock {\em Journal of Algebra}, 508, 05 2018.

\bibitem{boltje-pfeiffer}
R.~Boltje and G.~Pfeiffer.
\newblock An algorithm for the unit group of the {B}urnside ring of a finite
  group.
\newblock In {\em Groups St Andrews 2005}, volume 339, pages 230--236. London
  Math. Soc. Lectures Notes Series, 2007.

\bibitem{handbook}
S.~Bouc.
\newblock Burnside rings.
\newblock In {\em Handbook of {A}lgebra}, volume~2, chapter~6E, pages 739--804.
  Elsevier, 2000.

\bibitem{rittersegal}
S.~Bouc.
\newblock A remark on a theorem of {R}itter and {S}egal.
\newblock {\em J. of Group Theory}, 4:11--18, 2001.

\bibitem{dadegroup}
S.~Bouc.
\newblock The {D}ade group of a $p$-group.
\newblock {\em Inv. Math.}, 164:189--231, 2006.

\bibitem{burnsideunits}
S.~Bouc.
\newblock The functor of units of {B}urnside rings for $p$-groups.
\newblock {\em Comm. Math. Helv.}, 82:583--615, 2007.

\bibitem{bisetfunctorsMSC}
S.~Bouc.
\newblock {\em Biset functors for finite groups}, volume 1990 of {\em Lecture
  Notes in Mathematics}.
\newblock Springer-Verlag, Berlin, 2010.

\bibitem{both3}
S.~Bouc and J.~Th\'evenaz.
\newblock A sectional characterization of the {D}ade group.
\newblock {\em Journal of Group Theory}, 11(2):155--298, 2008.

\bibitem{boya}
S.~Bouc and E.~Yal{\c c}{\i}n.
\newblock Borel-{S}mith functions and the {D}ade group.
\newblock {\em J. Algebra}, 311:821--839, 2007.

\bibitem{dressresoluble}
A.~Dress.
\newblock A characterization of solvable groups.
\newblock {\em Math. Zeit.}, 110:213--217, 1969.

\bibitem{feit-thompson}
W.~Feit and J.~G. Thompson.
\newblock Solvability of groups of odd order.
\newblock {\em Pacific J. Math.}, 13:775--1029, 1963.

\bibitem{gluck}
D.~Gluck.
\newblock Idempotent formula for the {B}urnside ring with applications to the
  $p$-subgroup simplicial complex.
\newblock {\em Illinois J. Math.}, 25:63--67, 1981.

\bibitem{matsuda}
T.~Matsuda.
\newblock On the unit group of {B}urnside rings.
\newblock {\em Japan. J. Math.}, 8(1):71--93, 1982.

\bibitem{matsudamiyata}
T.~Matsuda and T.~Miyata.
\newblock On the unit groups of the {B}urnside rings of finite groups.
\newblock {\em J. Math. Soc. Japan}, 35(1):345--354, 1983.

\bibitem{ritter}
J.~Ritter.
\newblock Ein {I}nduktionssatz f\"ur rational {C}haraktere von nilpotenten
  {G}ruppen.
\newblock {\em J. f. reine u. angew. Math.}, 254:133--151, 1972.

\bibitem{segal}
G.~Segal.
\newblock Permutation representations of finite $p$-groups.
\newblock {\em Quart. J. Math. Oxford}, 23:375--381, 1972.

\bibitem{tomdieckgroups}
T.~tom Dieck.
\newblock {\em Transformation groups and representation theory}, volume 766 of
  {\em Lecture Notes in Mathematics}.
\newblock Springer-Verlag, 1979.

\bibitem{yalcin}
E.~Yal{\c c}{\i}n.
\newblock An induction theorem for the unit groups of {B}urnside rings of
  2-groups.
\newblock {\em J. Algebra}, 289:105--127, 2005.

\bibitem{yoshidaidemp}
T.~Yoshida.
\newblock Idempotents of {B}urnside rings and {D}ress induction theorem.
\newblock {\em J. Algebra}, 80:90--105, 1983.

\bibitem{yoshidaunit}
T.~Yoshida.
\newblock On the unit groups of {B}urnside rings.
\newblock {\em J. Math. Soc. Japan}, 42(1):31--64, 1990.

\end{thebibliography}

\vspace{3ex}
\begin{flushleft}
Serge Bouc\\
LAMFA-CNRS UMR 7352\\
Universit\'e de Picardie-Jules Verne\\
33, rue St Leu, 80039 Amiens Cedex 01\\
France\\
{\small\tt serge.bouc@u-picardie.fr}\\
{\small\tt http://www.lamfa.u-picardie.fr/bouc/}
\end{flushleft}
\end{document}